\documentclass[12pt,a4paper]{amsart}
\usepackage{amsmath,amssymb}

 \usepackage[dvips]{color,graphicx}

\newcommand{\beqn}{\begin{eqnarray}}
\newcommand{\eeqn}{\end{eqnarray}}
\newcommand{\ep}{\varepsilon}

\newcommand{\la}{\lambda}

\newcommand{\C}{\mathbb{C}}

\newcommand{\Z}{\mathbb{Z}}

\newcommand{\F}{\mathcal{F}}

\newcommand{\B}{\mathbb{B}}
\DeclareMathOperator{\Pf}{Pf}

\theoremstyle{plain}
\newtheorem{thm}{Theorem}[section]
\newtheorem{lem}[thm]{Lemma}
\newtheorem{prop}[thm]{Proposition}

\newtheorem{Def}[thm]{Definition}
\newtheorem{exam}[thm]{Example}

\def\maru#1{{\ooalign{\hfill$\scriptstyle#1$\hfill\crcr$\bigcirc$}}}

\title[mixed expansion formula for the rectangular Schur functions]{Mixed expansion formula for the rectangular Schur functions
and the affine Lie algebra $A_1^{(1)}$}
\author[Ikeda, Mizukawa, Nakajima and Yamada]{Takeshi Ikeda, Hiroshi Mizukawa, Tatsuhiro Nakajima\\
and \\ Hiro-Fumi Yamada
}
\address{Takeshi Ikeda, Department of Applied Mathematics, Okayama University
of Science, Okayama 700-0005, 
Japan}
\email{ike@xmath.ous.ac.jp}
\address{Hiroshi Mizukawa, Department of Mathematics,  National Defense
Academy,  Yokosuka 239-8686, Japan}
\email{mzh@nda.ac.jp}
\address{Tatsuhiro Nakajima, Faculty of Economics, Meikai University, Urayasu 279-8550, 
Japan}
\email{tatsu.nkjm@mac.com}
\address{Hiro-Fumi Yamada, Department of Mathematics, Okayama University, Okayama 700-8530, 
Japan}
\email{yamada@math.okayama-u.ac.jp}
\date{}

\begin{document}
\begin{abstract}
Formulas are obtained
that express the Schur $S$-functions indexed by Young diagrams
of rectangular shape
as linear combinations of 
``mixed"
products 
of Schur's $S$- and $Q$-functions.
The proof is achieved by using 
representations of the affine Lie algebra of type $A_1^{(1)}.$
A realization of the basic representation
that is of ``$D_2^{(2)}$''-type 
plays the central role.
\end{abstract}
 \maketitle
\section{Introduction}
We derive combinatorial formulas  that express the Schur
$S$-functions indexed by Young diagrams of rectangular shape, the rectangular
$S$-functions for short, as linear combinations of 
``mixed"
products of 
$S$- and $Q$-functions.

The rectangular $S$-functions  
are studied in \cite{iy,my} 
from a viewpoint of representations of
the affine Lie algebra of type $A_1^{(1)}$ and $A_2^{(2)}$.
These functions appear 
as certain distinguished weight vectors in the so called \textit{homogeneous} 
realization of the basic representation $L(\Lambda_0)$
of $A^{(1)}_{1}$ (see \cite{kac}).
On the other hand, the Schur $Q$-functions
arise naturally in the representation of $D_{\ell+1}^{(2)}$-type
Lie algebras (\cite{ny}).
In the subsequent pursuit of various realizations of $L(\Lambda_0)$, 
our formula has come out 
as an application of 
the isomorphism $D_2^{(2)}\cong A_1^{(1)}.$
Roughly speaking, we can realize 
the space $L(\Lambda_0)$ as a tensor product of the spaces of
the Schur $S$- and $Q$-functions.
We call such a ``mixed'' realization as the homogeneous realization of type $D_2^{(2)}.$

Let us describe our main result 
in more detail.
Let $\mu$ be a partition and
$S_\mu(t)$ be the corresponding Schur $S$-function,
where $t=(t_1,t_2,t_3,\ldots)$, 
and 
each $t_j\;(j=1,2,\ldots)$ is the $j$-th 
\textit{power sum} $p_j$ divided by $j.$
Let $Q_{\la}(t)$ denote the Schur $Q$-function
indexed by a \textit{strict} partition $\la$, where
$t=(t_1,t_3,t_5,\ldots)$.
Let $\square(m,n)$ denote the Young diagram of the rectangular shape $(n^m).$
Set also $S_\mu(t^{(2)})=S_\mu(t_2,t_4,t_6,\ldots).$ 
Note that the set $$\{Q_\la(t)S_\mu(t^{(2)});{\la\ \text{is\ a\ strict\ partition\ and\ }\mu\text{\ is\ a\ partition}}\}$$
forms a basis of the space of the symmetric functions.

Let $m,n$ be non-negative integers.
Our formula (Theorem \ref{main}), called "mixed expansion formula", reads:
\begin{align}\label{introformula}
\sum_{\mu}\delta(\mu) Q_{\mu{[0]}}(t)S_{\mu{[1]}}(t^{(2)})=S_{\square(m,n)}(t),
\end{align}
where the summation runs over
a certain finite set of strict partitions
determined by $m$ and $n$.  
For each strict partition $\mu$, one associates a strict partition $\mu[0]$,
 a partition $\mu[1]$ and
a sign $\delta(\mu) = \pm 1$  in a combinatorial way.
We prove the  formula (\ref{introformula}) by 
comparing two realizations
of $L(\Lambda_0)$ mentioned above.
The left hand side stems from combinatorial
descriptions of actions of Chevalley generators 
in the homogeneous realization of type $D_2^{(2)},$
whereas the right hand side 
is obtained via ``vertex operator calculus'' (as employed in \cite{iy})
in the homogeneous realization of type $A_1^{(1)}.$

Here we explain the background of our study of rectangular Schur functions.
As written in the above,  our formula arose from a study of the homogeneous realization
of the basic $A^{(1)}_1$-module.  We have two pictures of the
principal realization of the basic  $A^{(1)}_1$-module; one is described
in terms of the 2-reduced Schur functions and is relevant to the KdV
hierarchy; the other is the twisted version, which is best described
by the $Q$-functions and is relevant to the 4-reduced BKP hierarchy.
On the other hand, the homogeneous realization of that module is
connected with the nonlinear Schr\"{o}dinger (NLS) hierarchy.
Using an intertwining operator between the (non twisted) principal
and the homogeneous realizations, one can derive an expression of
the rectangular Schur functions and certain $\tau$-functions of the
NLS hierarchy (\cite{iy}). 

     The paper is organized as follows.
     In Section \ref{partition} we recall some combinatorial materials related to partitions. 
     In Section \ref{formula} we state our main theorem on rectangular Schur functions.
     In Section \ref{spin} we recall the spin representation of 
     $A_{1}^{(1)}$
    and describe the action of $A_{1}^{(1)}$  in terms of Young diagrams.
In Section \ref{bozon} through the boson-fermion correspondence, we
obtain weight vectors as a sum of products of $S$- and $Q$- functions.
In Section \ref{vertex} we consider $f_{i}$-action ($i=0,1$) and obtain the
rectangular Schur functions  appearing in the right hand side of our formula through a vertex operator calculus.  
Section \ref{proof} is devoted to the proof of the main theorem.

\section{Combinatorics of Partitions}\label{partition}
\subsection{Partition}
A \textit{partition} is any non-increasing sequence of non-negative
integers $\lambda=(\lambda_1,\lambda_2,\ldots)$
containing only finitely many non-zero terms.
We regard two partitions as the same that differ only by a string of
zeros at the end.
The non-zero $\lambda_i$ are called the \textit{parts}
  of $\lambda.$
The number of parts is the \textit{length} of $\lambda$, denoted by  
$\ell(\lambda).$
A partition is \textit{strict} if all parts are distinct.
Denote by $\mathcal{P}$ (resp. $\mathcal{SP}$)
the set of all partitions (resp. strict partitions).

\subsection{$i$-addable node}
Let $\lambda\in \mathcal{SP}.$
To each node $x \in \lambda$ in the $j$-th column,
we assign a color $c(x)$ by the following rule:
$$
c(x)=\begin{cases}0\quad (j\equiv 0,1 \mod 4)\\
1 \quad (j\equiv 2,3 \mod 4)
\end{cases}.
$$
For example, the nodes of $\lambda=(5,4,2,1)$ is colored as
$$\begin{array}{ccccc}0&1&1&0&0\\
0&1&1&0&\\
0&1&&&\\
0&&&&
\end{array}.
$$ 
We say that a node $x$ 
is $i$-addable to $\lambda$,
if $\lambda\cup \{x\}$ is a strict partition and $c(x)=i.$
The following nodes indicated by dots are the $1$-addable nodes:
$$\begin{array}{cccccc}0&1&1&0&0&\bullet\\
0&1&1&0&&\\
0&1&\bullet&&&\\
0&&&&
\end{array}.
$$
Set
$${I}_i^\ell(\lambda)
=\{\mu\in \mathcal{SP}\;|\;
\mu\supset \lambda,\;|\mu|=|\lambda|+\ell,\;\forall x\in \mu-\lambda,  
\;c(x)=i\}.
$$
It is the set of strict partitions
obtained from $\lambda$
by adding $i$-nodes $\ell$ times in succession.
Put
$$
\begin{cases}
c_m=(4m-3,\ldots,5,1) & (m>0)\\
c_{m}=\emptyset &  (m=0)\\
c_m=(-4m-1,\ldots,7,3) & (m<0).
\end{cases}
$$
If $m>0$, we have ${I}^\ell_1(c_m)=\emptyset$ for $\ell > 2m$
and ${I}^\ell_0(c_{-m})=\emptyset$ for $\ell > 2m+1$.

The strict partitions $c_m\;(m\in \Z)$ are called
{\it $4$-bar cores},
introduced in \cite{bol, ny}. 

\begin{exam}
For $m=-2$ and $i=0$ we have
$$I_{0}^{1}(c_{-2})=\{(8,3),(7,4),(7,3,1)\},$$
$$
{I}_0^2(c_{-2})=
\{
(9,3),(8,4),(8,3,1),(7,4,1),(7,5)
\}
$$
and
$$
{I}_0^3(c_{-2})=
\{
(9,4),(8,5),(9,3,1),(8,4,1),(7,5,1)
\}.
$$
\end{exam}

\subsection{$4$-bar quotient}
Let us introduce the notion of $4$-bar quotient.
We shall give a bijection
$$\mathcal{SP}\rightarrow\Z\times \mathcal{SP}\times \mathcal{P},\quad
\lambda\mapsto (m,\lambda[0],\lambda[1]).$$
For $\lambda\in \mathcal{SP}$,
the pair $(\lambda[0],\lambda[1])$ is called the $4$-bar quotient
of $\lambda.$

Let us identify the strict partition $\lambda$ with the subset
$\pmb{\lambda}=\{\lambda_1,\ldots,\lambda_s\}$ of $\mathbb{N}.$
For $a=0,1,2,3$, we set $\pmb{\lambda}^{(a)}=\{\lambda_j\in
\pmb{\lambda}\,|\,\lambda_j\equiv a \mod 4\}.$
Namely
$$\pmb{\lambda}^{(a)}=\pmb{\lambda}\cap(4\mathbb{N}+a)\quad(a=0,1,2,3)$$
and
we have $\pmb{\lambda}=\sqcup_{a=0}^3\,{\pmb{\lambda}}^{(a)}.$
The even part $\pmb{\lambda}^{(0)}\cup\pmb{\lambda}^{(2)}
\subset 2\mathbb{N}$ of $\pmb{\lambda}$
gives a strict partition $\lambda[0]$
via the inclusion
$$
\pmb{\lambda}^{(0)}\cup\pmb{\lambda}^{(2)}\subset 2\mathbb{N}
\longrightarrow\mathbb{N},\quad
2k\mapsto k.
$$
From the odd parts $\pmb{\lambda}^{(1)},\pmb{\lambda}^{(3)}$,
we define a partition $\lambda[1]$ in the following way:
First consider two bijections
$$
\iota:4\mathbb{N}+1\longrightarrow\Z_{\geq 0}\quad(4k+1\mapsto k),\quad
\iota^*:4\mathbb{N}+3\longrightarrow\Z_{< 0}\quad(4k+3\mapsto -k-1).
$$
Then define a subset
$$
\mathcal{M}(\lambda)=\iota(\pmb{\lambda}^{(1)})
\cup(\Z_{<0}- \iota^*(\pmb{\lambda}^{(3)}))
$$
of $\Z.$ This is a ``Maya diagram'' in the sense that, if we
express $\mathcal{M}(\lambda)$ as an descending sequence
$i_1>i_2>i_3>\cdots,$
then the integer $m=\sharp \pmb{\lambda}^{(1)}-
\sharp \pmb{\lambda}^{(3)}$ satisfies $i_k=-k+m$ for $k\ll 0.$
Then we can define
$$
\lambda[1]=(i_1+1-m,i_2+2-m,i_3+3-m,\ldots)\in\mathcal{P}.
$$
The integer $m$ is called the {\it charge} of $\mathcal{M}(\lambda).$
\begin{lem}\label{SPZSPP} (cf. \cite{bol})
The map
$$\mathcal{SP}\rightarrow\Z\times \mathcal{SP}\times \mathcal{P},\quad
\lambda\mapsto (m,\lambda[0],\lambda[1])$$
is a bijection.
\end{lem}
       %

We can illustrate the above construction.
Let us look at a particular example, $\lambda=(11,9,6,2,1,0)$.
We draw a ``4-bar abacus":
\begin{align*}
 {\begin{array}{ccc}
0&\maru{1}&3\\
\maru{2}&&\\
4&5&7\\
\maru{6}&&\\
8&\maru{9}&\maru{11}\\
10&&\\
12&13&15\\
\end{array}}
\end{align*}
Here we do not put a bead on $0$.
We can read ${\pmb{\lambda}}^{(0)}\bigcup {\pmb{\lambda}}^{(1)}$ from the
first column. Then we have
$\lambda[0]=(3,1).$
From the second and the  third columns, we can read
\begin{align*}
\iota({\pmb \lambda}^{(1)})&=(2,0)\\
\iota^*({\pmb \lambda}^{(3)})&=(-3).\end{align*}
Then we obtain 
$$M=(2,0,-1,-2,-4,-5,\cdots)$$
and draw a Maya diagram;
\begin{figure}[h]
\centering
\unitlength=0.1mm
\begin{picture}(500,300)
\thicklines
\put(-30,-30){\line(1,1){230}}
\put(80,-0){\line(1,1){160}}
\put(240,80){\line(1,1){40}}
\put(80,0){\line(-1,1){40}}
\put(120,40){\line(-1,1){40}}
\put(160,80){\line(-1,1){40}}
\put(240,80){\line(-1,1){80}}
\put(400,0){\line(-1,1){200}}
%
\put(200,200){\line(0,1){40}}
\put(240,200){\line(0,1){40}}
\put(280,200){\line(0,1){40}}
\put(320,200){\line(0,1){40}}
\put(360,200){\line(0,1){40}}
\put(400,200){\line(0,1){40}}
\put(0,200){\line(0,1){40}}
\put(40,200){\line(0,1){40}}
\put(80,200){\line(0,1){40}}
\put(120,200){\line(0,1){40}}
\put(160,200){\line(0,1){40}}
\put(-40,200){\line(1,0){460}}
\put(-40,240){\line(1,0){460}}
\put(-37,202){\LARGE{$\bullet$}}
\put(3,202){\LARGE{$\bullet$}}
\put(83,202){\LARGE{$\bullet$}}
\put(123,202){\LARGE{$\bullet$}}
\put(163,202){\LARGE{$\bullet$}}
\put(243,202){\LARGE{$\bullet$}}
\put(-38,246){\tiny{$-\!5$}}
\put(1,246){\tiny{$-\!4$}}
\put(81,246){\tiny{$-\!2$}}
\put(121,246){\tiny{$-\!1$}}
\put(173,246){\tiny{$0$}}
\put(253,246){\tiny{$2$}}
\end{picture}.
\end{figure}

\noindent
Finally  we have
$\lambda[1]=(2,1,1,1)$
and $m=1$.

\subsection{Sign}
Each strict partition $\mu$ in ${I}^{\ell}_1(c_{m})$
or ${I}^{\ell}_0(c_{m})$ has its own sign determined 
by bead configuration. 
\begin{Def}\label{gabacus}
Put a bead on $0$
of the 4-bar abacus of $\lambda \in {I}^{\ell}_i(c_{m})\ (i=0,1)$,
 if and only if $m<0$ and $\lambda_{-m+1}=0$.
Let $g(\lambda)$ be  the number of  pair of
 beads on the central runner at the positions bigger than that
of each bead on the leftmost runner.
For a strict partition $\lambda \in {I}^{\ell}_i(c_{m})$, we define 
the sign  by
$$\delta(\lambda)=(-1)^{g(\lambda)}.$$ 
\end{Def}
\begin{exam}
We consider the case of $i=1, m=3,\ell=3$ and $\lambda=(11,5,2)$.
\begin{align*}
 {\begin{array}{ccc}
0&1&3\\
\maru{2}&&\\
4&\maru{5}&7\\
6&&\\
8&9&\maru{11}\\
\end{array}}
\end{align*}
We have $\delta(\lambda)=(-1)^{1}=-1$.
In the case of $i=0,m=-4,\ell=5$ and $\lambda=(15,13,8,5)$,
we have to put a bead on 0.
\begin{align*}
 {\begin{array}{ccc}
\maru{0}&1&3\\
2&&\\
4&\maru{5}&7\\
6&&\\
\maru{8}&9&11\\
10&&\\
12&\maru{13}&\maru{15}\\
\end{array}}
\end{align*}
We have $\delta(\lambda)=(-1)^{2+1}=-1$.
\end{exam}


\section{Main Result}\label{formula}

Define $h_n(t)$ by
$
\exp\left(\sum_{n=1}^\infty t_{n}z^{n}\right)
=\sum_{n=0}^\infty h_n(t)z^n.
$
Let $\lambda$ be a partition.
The Schur $S$-function with shape $\lambda$
is defined as
$$S_{\lambda}(t)=\det(h_{\lambda_i+j-i}(t)).
$$
Define $q_n(t)$ by $\exp\left(\sum_{n=1}^\infty t_{2n-1}z^{2n-1}\right)
=\sum_{n=0}^\infty q_n(t)z^n.$
For
$m>n\geq 0,$ we put
$$Q_{m,n}(t)=q_m(t)q_n(t)+2\sum_{i=1}^n(-1)^iq_{m+i}(t)q_{n-i}(t).$$
If $m\leq n$ we define $Q_{m,n}(t)=-Q_{n,m}(t).$
Let $\lambda=(\lambda_1,\ldots,\lambda_{2n})$ be a strict partition,
where $\lambda_1>\cdots>\lambda_{2n}\geq 0.$
Then the
$2n\times 2n$ matrix $M_\lambda=(Q_{\lambda_i,\lambda_j})$ is
skew-symmetric.
The $Q$-function $Q_\lambda$ is defined as
$$
Q_\lambda(t)=\mathrm{Pf}(M_\lambda).
$$

We can now state our main result which we call the {\it mixed expansion formulas}.

\begin{thm}\label{main}For non-negative integers $m$ and $n$, we have

$$\sum_{\mu \in {I}_1^n(c_m)}
\delta(\mu) Q_{\mu[0]}(t) S_{\mu[1]}(t^{(2)})
= S_{\square (2m-n,n)} (t), $$

$$\sum_{\mu \in {I}_0^n(c_{-m})}
\delta(\mu) Q_{\mu[0]}(t) S_{\mu[1]}(t^{(2)})
= S_{\square (n,2m+1-n)} (t),$$
where $t=(t_1,t_2,t_3,\cdots)$ and $S_{\nu}(t^{(2)})=S_{\nu}(u)|_{u_j \mapsto t_{2j}}.$
If ${I}_i^n(c_{\pm m})=\emptyset$, we agree that  both sides are  equal to 0. 
\end{thm}
\begin{exam}
For $m=3$ and $n=2$, we have 
$$
S_{1^4}(t^{(2)})+S_{21^2}(t^{(2)})+S_{2^2}(t^{(2)})+Q_{53}(t)-Q_{51}(t)S_{1}(t^{(2)})+Q_{31}(t)S_{2}(t^{(2)})
= S_{\square (4,2)} (t).$$
For $m=-2$ and $n=2$, we have 
$$
-S_{3}(t^{(2)})-S_{21}(t^{(2)})+S_{1^2}(t^{(2)})Q_{2}(t)+S_{1}(t^{(2)})Q_{4}(t)+Q_{42}(t)
= S_{\square (2,3)} (t).$$
\end{exam}

\section{The spin representation of $A_1^{(1)}$}\label{spin}

We consider the associative $\C$-algebra $\B$
defined by the generators $\beta_n\;(n\in \Z)$
and the anti-commutation relations:
$$[\beta_m,\beta_n]_{+}=\beta_m\beta_n+\beta_n\beta_m=(-1)^m\delta_{m+n,0}.$$

These generators are often called the
neutral free fermions.  
Note that $\beta_0^2=1/2$. 
Let $\mathcal{F}$ be the Fock module which is a left
$\B$-module generated by  the vacuum $| \emptyset \rangle$ with 
$$\beta_n|\emptyset\rangle=0 \quad (n<0).$$
Similarly we consider the right $\B$-module
$
\F^\dagger
$
which is generated by the vacuum $\langle \emptyset|$ with
$$
\langle \emptyset|\beta_n=0 \quad (n>0).
$$
Elements of $\mathcal{F}$ and $\mathcal{F}^{\dagger}$ are sometimes called ``states".
We have a bilinear pairing
$$\F^\dagger\otimes_\B \F\rightarrow \C,\quad
\langle \emptyset | u\otimes_\B v|\emptyset\rangle\mapsto \langle \emptyset| uv |\emptyset\rangle.$$
This pairing is called the {\it vacuum expectation value}. 
The vacuum expectation value  is uniquely determined by putting
$\langle\emptyset|\emptyset\rangle=1$ and $ \langle \emptyset|\beta_0|\emptyset\rangle=0$.


\begin{Def}
Let $\lambda$ be a strict partition, which we may write
in the form $\lambda=(\lambda_1,\ldots,\lambda_{2k})$ where
$\lambda_1>\cdots>\lambda_{2k}\geq 0.$
Here 
\begin{align*}2k=
\begin{cases} 
\ell(\lambda) & \ell(\lambda)\equiv 0\pmod{2},\\
\ell(\lambda) +1 & \ell(\lambda)\equiv 1\pmod{2}.
\end{cases}
\end{align*}
Let $|\lambda\rangle$ denote the state
$$
|\lambda\rangle=\beta_{\lambda_1}\cdots\beta_{\lambda_{2k}}
|\emptyset\rangle \in \F.
$$
For $\lambda=\emptyset$, we define $|\emptyset \rangle=|\emptyset\rangle$.
\end{Def}

Set $f^\infty_i=(-1)^i\beta_{i+1}\beta_{-i}\;(i\geq 0).$ They have
the following combinatorial property whose proof is left to the reader.
\begin{prop} Let $\lambda$ be a strict partition. 
 If $i > 0$, then we have
     $$
     f^\infty_i | \lambda\rangle
       =
     \begin{cases}
        | \mu \rangle &\text{if $i$ is a part of $\lambda$ and $i+1$ is
not},
      \\
        0 &\text{otherwise},
     \end{cases}
$$
where $\mu$ is obtained from $\lambda$ by replacing
its part $i$ by $i+1$.
If $1$ is not a part of $\lambda$, then
we have
$$f^\infty_0|\lambda\rangle
=\begin{cases}{2}^{-1}|\mu\rangle&
\text{$1$ is not a part of $\lambda$ and\ }
\ell(\lambda)\equiv 1\pmod{2},\\
|\mu\rangle &
\text{$1$ is not a part of $\lambda$ and\ }
 \ell(\lambda)\equiv 0\pmod{2},\\
 0 & otherwise,
\end{cases}$$
where 
$\mu$ is obtained from $\lambda$ by adding
a part $1$.
\end{prop}

We shall use standard notation
of the affine Lie algebra $A_1^{(1)} $ (\cite{kac}).
Let $e_i,f_i,h_i(i=0,1)$ be the Chevalley generators,
$\alpha_0,\alpha_1$ are the simple roots,
 $\delta=\alpha_{0}+\alpha_{1}$ is the fundamental imaginary root,
 $\Lambda_i(i=0,1)$ are the fundamental weights.
The affine Lie algebra $A^{(1)}_1$ acts on  $\F$
by
$$
f_0=\sqrt{2}\sum_{n\in \Z}\beta_{-4n+1}\beta_{4n},\quad
f_1=\sqrt{2}\sum_{n\in \Z}\beta_{-4n-1}\beta_{4n+2},
$$
$$
e_0=\sqrt{2}\sum_{n\in \Z}\beta_{4n}\beta_{-4n-1},\quad
e_1=\sqrt{2}\sum_{n\in \Z}\beta_{4n-2}\beta_{-4n+1},
$$
$$
h_1=-h_0+1=2\sum_{n\in \Z}:\beta_{4n-1}\beta_{-4n+1}:,
$$
where we define the normal ordering for the quadratic
elements by 
$$:\beta_n\beta_m:=\beta_n\beta_m-\langle\emptyset|\beta_n\beta_m|\emptyset \rangle.$$
Let $\F_{0}$ 
 be the $A_{1}^{(1)}$-submodule of $\F$
generated by $|\emptyset\rangle$
.
$\F_0$ is isomorphic to the irreducible highest weight module
$L(\Lambda_0).$

Note  the following expressions:
$$f_0=\sqrt{2}\sum_{j\geq 0}f^\infty_{4j}+\sqrt{2}\sum_{j\geq  
0}f^\infty_{4j+3},\quad
f_1=\sqrt{2}\sum_{j\geq 0}f^\infty_{4j+1}+\sqrt{2}\sum_{j\geq  
0}f^\infty_{4j+2}.
$$

We need the following combinatorial lemmas:
\begin{lem}\label{hardcore}\cite{ny} A weight vector of the weight  
$\Lambda_0-m^2\delta+m\alpha_1$ is
given by $|c_m\rangle$ in ${\F_{0}}.$
\end{lem}
The weight diagram of $L(\Lambda_0)$ looks as follows.  Maximal weights correspond
to the lattice points on the parabola, and other weights are on the lattice points under this parabola.
\\

\begin{center}{
\includegraphics[width=6cm]{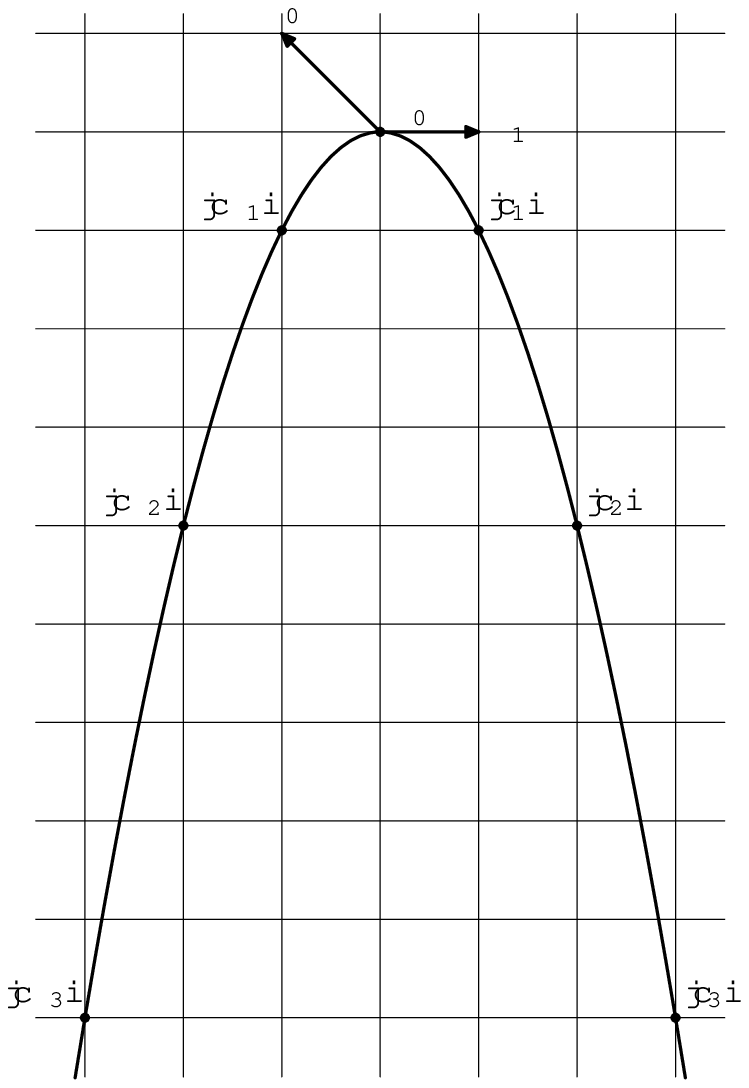}
}\end{center}

\begin{lem}\label{co}
$$
\frac{f_i^\ell}{\ell!}\,|c_m\rangle
=\sqrt{2}^{-\varepsilon_m}
\sum_{\lambda \in {I}_i^\ell(c_m)}
\sqrt{2}^{\,a(\lambda)}|\lambda\rangle.
$$
where $
a(\lambda)=\sharp\{j\;|\;\lambda_j\equiv 0\mod 2\}
$
for
 $\lambda=(\lambda_{1},\lambda_{2},\cdots,\lambda_{\ell(\lambda)+\ep_{\ell(\lambda)}})$
with $\ep_m=1$ if $m$ is odd, and $\ep_m=0$ if $m$ is even.
\end{lem}

\begin{proof}
Firstly we consider the case of $i=1$.
For $\lambda \in {I}_{1}^{\ell}(c_m)$ $(m>0)$, put 
\begin{align*}\lambda - c_{m}=
 (r_{1},r_{2},\cdots,r_{m}).
\end{align*}
We compute
 \begin{align*}
 r_{1}!r_{2}!\cdots r_{m}!&=2^{\sharp \left\{ j ; r_{j}=2\right\}}\\
 &=2^{\left(\ell-\sharp \left\{ j ; r_{j}=1\right\}\right)/2}\\
 &=2^{(\ell-a(\lambda)+\varepsilon_{m})/2},
 \end{align*}
where we note that 
$a(\lambda)$ counts a 0 if $m \equiv 1 \pmod 2$. 
Then  the coefficient of $|\lambda \rangle$ is
 \begin{align*}
 \frac{\sqrt{2}^{\ell}}{\ell !}  \frac{\ell !}{r_{1}!r_{2}!\cdots r_{m}!}
 =
 \sqrt{2}^{a(\lambda)-\varepsilon_m}.
   \end{align*}
 Secondly we consider the case of $i=0$.
 In this case we have to take $\beta_{1}\beta_{0}$-part into account.
 Let $$\lambda- c_{-m}=(r_{1},r_{2},\cdots,r_{m},r_{m+1})$$
 for $\lambda \in {I}_{0}^{\ell}(c_{-m})$ $(m \geq 0)$.
 Then the coefficient of $|\lambda \rangle$ is
 \begin{align}\label{coeff}
\frac{1}{2^{r_{m+1}\varepsilon_{m}}} \frac{\sqrt{2}^{\ell}}{\ell !}  \frac{\ell !}{r_{1}!r_{2}!\cdots r_{m}!}.
 \end{align}
A computation similar to the case $i=1$ above 
$$r_{1}!r_{2}!\cdots r_{m}!=2^{\sharp \left\{ j ; r_{j}=2\right\}}
=2^{\left(\ell-\sharp \left\{ j ; r_{j}=1, j \leq m \right\}-r_{m+1}\right)/2}.$$
Here we divide our argument into two cases.   
First we assume that $m$ is even. We have
 \begin{align*}
\sharp \left\{ j ; r_{j}=1, j \leq m \right\}&=
\begin{cases}
a(\lambda) & (\lambda_{m+1}=0)\\
a(\lambda)-1& (\lambda_{m+1}=1)
\end{cases}\\
 \end{align*}
Secondly we assume that $m$ is odd. We have
 \begin{align*}
\sharp \left\{ j ; r_{j}=1, j \leq m \right\}&=
\begin{cases}
a(\lambda)-1 & (\lambda_{m+1}=0)\\
a(\lambda)& (\lambda_{m+1}=1)
\end{cases}\\
 \end{align*}
By substituting these four results into (\ref{coeff}) we obtain 
 \begin{align*}
 \frac{1}{2^{r_{m+1}\varepsilon_{m}}}\frac{\sqrt{2}^{\ell}}{\ell !}  \frac{\ell !}{r_{1}!r_{2}!\cdots r_{m}!}
 =
 \sqrt{2}^{a(\lambda)-\varepsilon_m}.
   \end{align*}
\end{proof}

\begin{exam}
$
f_0|c_{-2}\rangle=\sqrt{2}\,\beta_8\beta_3|\emptyset\rangle
+\sqrt{2}\,\beta_7\beta_4|\emptyset\rangle
+\sqrt{2}\,\beta_7\beta_3\beta_{1}\beta_0|\emptyset\rangle.
$

\end{exam}

\section{Bosonization}\label{bozon}

In this section we will establish the bozon-fermion correspondence and see
the states as the polynomials. In the course, products of Schur's $S$- and $Q$-functions
arise naturally.    

We introduce the operators $\phi_n,\psi_n,\psi^*_n\;(n\in\Z)$ by
\begin{align}\label{kakikae}
\beta_{4n}=\phi_{2n},\quad
\beta_{4n+1}=\sqrt{-1}\,\psi_n,\quad
\beta_{4n+2}=\sqrt{-1}\,\phi_{2n+1},\quad
\beta_{4n+3}=\sqrt{-1}\,\psi^*_{-n-1},
\end{align}
which satisfy the anti-commutation
relations:
\begin{align*}
[\psi_m,\psi^*_n]_{+}&=\delta_{m,n},\quad
[\psi^*_m,\psi^*_n]_{+}=[\psi_m,\psi_n]_{+}=0,\\
[\phi_m,\phi_n]_{+}&=(-1)^m\delta_{m+n,0},\\
[\psi^*_m,\phi_n]_{+}&=[\psi_m,\phi_n]_{+}=0.
\end{align*}

Let us introduce the
bosonic current operators
$$
H_{2m}=\sum_{k\in\Z}\;:\psi_{k}\psi^*_{k+m}:\;,\quad
H_{2m+1}=\frac{1}{2}\,\sum_{k\in\Z}\;(-1)^{k+1}\phi_{k}\phi_{-k-(2m+1)}.
$$
These operators generate an infinite-dimensional Heisenberg algebra
 $$\mathfrak{H}=\oplus_{n\ne 0}\C H_n\oplus \C c,$$
where $c$ denotes the central element of $\mathfrak{H}$.
One has
$$
[H_m,H_n]=\frac{m}{2}\,\delta_{m+n,0}\,c.
$$

We have a canonical $\mathfrak{H}$-module
$S[\mathfrak{H}_{-}]$,
where $\mathfrak{H}_{-}=\oplus_{n<0}\C H_n,$ and
$S$ stands for the symmetric algebra.
Let $t_n=\frac{2}{n}H_{-n}\;(n>0).$ Then we can identify  
$S[\mathfrak{H}_{-}]$
with the ring $\C[\,t\,]=\C[\,t_1,t_2,t_3,\ldots\,]$ of
polynomials in infinitely many variables $t_n.$
The representation of $\mathfrak{H}$ on $\C[\,t\,]$ is described as  
follows:
$$
H_{n}P(t)=\frac{\partial}{\partial t_{n}}P(t),\quad
H_{-n}P(t)=\frac{n}{2}\,t_n P(t)\quad
(n>0,\;P(t)\in \C[\,t\,]),
$$
and $c$ acts as  identity.

If we introduce the space of highest weight vectors with respect to  
$\mathfrak{H}$
by $$
\Omega=\left\{|v\rangle \in\F\ ;\ H_m|v\rangle=0\;(\forall m>0)\right\},
$$
then $\Omega$ has a basis
$
\{|\sigma,m\rangle ;\ m\in \Z,\;\sigma=0,1\},
$
where
$$
|0,m\rangle=
\begin{cases}\psi_{m-1}\cdots\psi_0 |\emptyset\rangle
& (m>  0)\\
 |\emptyset\rangle
& (m= 0)\\
\psi^*_{m}\cdots\psi^*_{-1} |\emptyset\rangle
& (m<0)
\end{cases},\quad
|1,m\rangle=
\begin{cases}\sqrt{2}\phi_0\psi_{m-1}
\cdots\psi_0 |\emptyset\rangle
&(m> 0)\\
\sqrt{2}\phi_0|\emptyset\rangle
&(m= 0)\\
\sqrt{2}\phi_0\psi^*_{m}\cdots\psi^*_{-1} |\emptyset\rangle
&(m<0)
\end{cases}.
$$
Note that
$$
\phi_n|\sigma,m\rangle=0\quad(n<0),\quad
\psi_n|\sigma,m\rangle=0\quad(n<m),\quad
\psi^*_n|\sigma,m\rangle=0\quad(n\geq m).
$$

\begin{lem}\label{vacid}
$$|c_m\rangle=(\sqrt{-1})^{-|m|}\,\sqrt{2}^{- 
\varepsilon_m}|\,\varepsilon_m,m\rangle.$$
\end{lem}
\begin{proof}
We can easily obtain the equation by direct calculation.
\end{proof}

We introduce  formal symbols $\theta$ and  $e^{m\alpha}$ 
which satisfies $\theta^2=1$
and define
$$
\pmb{\Omega}=\bigoplus_{m\in \Z,\,\sigma=0,1}
\C\, \theta^\sigma e^{m\alpha}.
$$
Then $H_{n}$ act on $\C[\,t\,]\otimes \pmb{\Omega}$ by
$H_{n}\otimes id$.
\begin{prop} \cite{djkm1,djkm2} There exists a canonical isomorphism of  
$\mathfrak{H}$-modules
$$
\Phi:\F\longrightarrow
\C[\,t\,]\otimes_\C \pmb{\Omega}
$$
such that $
\Phi(|\sigma,m\rangle)=\theta^\sigma e^{m\alpha}(m\in \Z,\sigma=0,1).$
\end{prop}
We will see  $\C[\,t\,]\otimes \pmb{\Omega}$ as an $A_{1}^{(1)}$-module 
via $\Phi$ (cf. Proposition \ref{intertw}).

When we write
$\Phi(|v\rangle)=\sum_{m,\sigma}P_{m,\sigma}(t)\,\theta^\sigma  
e^{m\alpha}$
for $|v\rangle\in \F$,
the coefficient $P_{m,\sigma}(t)\in \C[t]$
can be expressed
in terms of the vacuum expectation value on $\B$
as follows:
$$
P_{m,\sigma}(t)=\langle m,\sigma|e^{H(t)}|v\rangle,\quad
H(t)=\sum_{n=1}^\infty t_{n}H_{n}.
$$
Introduce the states $\langle m,\sigma|\in \F^\dagger\;(m\in  
\Z,\sigma=0,1)$ which are
characterized by
$
\langle m,\sigma|\sigma',n \rangle=\delta_{m,n}\delta_{\sigma,\sigma'}\;
(m,n\in \Z,\sigma,\sigma'=0,1)
$
and
$$
\langle m,\sigma|\phi_n=0\quad(n>0),\quad
\langle m,\sigma|\psi_{n}=0\quad( n\geq m),\quad
\langle m,\sigma|\psi^*_{n}=0\quad(n<m).
$$

We denote by $\mathbb{W}_\phi$ the linear subspace of $\B$
spanned by $\phi_n\;(n\in \Z).$
\begin{lem}\label{PF}
If
$
\langle u|\in \F^\dagger,|v\rangle \in \F
$
be such that
$\langle u|\phi_{n}=0\;(n>0),\;
\phi_n|v\rangle =0\;(n<0)
$, then for $w_i\in \mathbb{W}_\phi\;(i=1,\ldots,2k)$  we have
$$
\langle u|w_{1}\cdots w_{2k}
|v\rangle
=\langle u|v\rangle
\Pf(\langle \emptyset| w_i w_j|\emptyset \rangle)
$$
\end{lem}
\begin{proof}
A bilinear form on $\mathbb{W}_\phi$ is
defined by $(a,b)\mapsto \langle u|ab|v\rangle$,
which has all the properties of vacuum expectation value on $\mathbb{W}_\phi$
except for the normalization condition.
Obviously, the normalization factor is given by $\langle u|v\rangle.$
Hence the lemma follows.
\end{proof}

\begin{lem}\label{BFC}\cite{djkm1,djkm2}
We have
\begin{align*}
&\Phi(\psi_{i_1}\cdots\psi_{i_s}|0,m\rangle)
=S_{(i_1-m,i_2-m,\ldots,i_s-m)-\delta_s}(t^{(2)})\, e^{(m+s)\alpha}
\quad(i_1>\cdots>i_s> m),\\
&\Phi(\phi_{j_1}\cdots\phi_{j_a}|\emptyset\rangle)=
\sqrt{2}^{-a}Q_{j_1,\ldots,j_a}(t)\,\theta^{a}
\quad(j_1>\cdots>j_a\geq 0),
\end{align*}
where $\delta_s=(s-1,s-2,\ldots,1,0)$
and $t^{(2)}=(t_2,t_4,\ldots).$
\end{lem}

Lemma \ref{PF} and \ref{BFC} give us

\begin{lem}\label{QS} Let $j_1>\cdots>j_a\geq 0,\;i_1>\cdots>i_s>m.$
We have
$$
\Phi(\phi_{j_1}\cdots \phi_{j_a}\psi_{i_1}\cdots\psi_{i_s}|0,m\rangle)
=\sqrt{2}^{-a}\,Q_{j_1,\ldots,j_a}(t)\,S_{(i_1-m,\ldots,i_s-m)- 
\delta_s}(t^{(2)})
\theta^{a} e^{(m+s)\alpha}.
$$
\end{lem}

Consequently we  obtain the following proposition.

\begin{prop}
Let $\lambda\in I^\ell_i(c_m).$ There exists a  $4$-th root of
unity $\zeta_{m,\ell,i}(\lambda)$ such that
$$
\Phi(\sqrt{2}^{\,a(\lambda)}|\lambda\rangle)=\zeta_{m,\ell,i}(\lambda)\;
\,
Q_{\lambda[0]}(t)S_{\lambda[1]}(t^{(2)})\,
\theta^{m+\ell}\, e^{(m+(-1)^i\ell)\alpha}.$$
\end{prop}

\section{Vertex operators}\label{vertex}

In this section we realize $f_i$ on $\B$ in terms of vertex operators.
We introduce the formal generating functions
$$
\phi(z)=\sum_{n\in \Z}\phi_{n} z^{n},\quad
\psi(z)=\sum_{n\in \Z}\psi_n z^{2n},\quad
\psi^*(z)=\sum_{n\in \Z}\psi^*_{-n} z^{2n-2}.
$$

For $t=(t_1,t_2,t_3,\ldots),$ set
$$\xi(t,z)=\sum_{n=1}^\infty t_nz^n,\quad
\xi_{0}(t,z)=\sum_{n=1}^\infty t_{2n}z^{2n},\quad
\xi_{1}(t,z)=\sum_{n=1}^\infty t_{2n-1}z^{2n-1}.
$$

On the space $\pmb{\Omega}$, we define the operators $\theta,e^{\pm  
\alpha}$ and $z^{H_{0}}$by
\begin{align*}
&\theta.(\theta\,e^{m\alpha})=e^{m\alpha},&
\theta.\,e^{m\alpha}=\theta \,e^{m\alpha},\quad\\
&e^{\pm\alpha}.(\theta\,e^{m\alpha})=-\theta\,e^{(m\pm 1)\alpha},&
e^{\pm\alpha}.e^{m\alpha}=e^{(m\pm 1)\alpha},\\
\end{align*}
and
$$z^{H_{0}}.(\theta\,e^{m\alpha})=z^m(\theta\,e^{m\alpha})\ \ (\sigma=0,1).$$


\begin{prop}\label{intertw}\cite{djkm1,djkm2}
One has
\begin{align*}
\Phi\; \phi(z)\Phi^{-1}&={\sqrt{2}^{-1}}e^{\xi_{1}(t,z)}
e^{-2\xi_{1}(\widetilde{\partial_t},z^{-1})}\theta,\\
\Phi\; \psi(z)  
\Phi^{-1}&=e^{\xi_{0}(t,z)}e^{-2\xi_{0}(\widetilde{\partial_{t}},z^{- 
1})}e^{\alpha}z^{2H_0},\\
\Phi\; \psi^*(z)  
\Phi^{-1}&=e^{-\xi_{0}(t,z)}e^{2\xi_{0}(\widetilde{\partial_{t}},z^{- 
1})}e^{-\alpha}z^{-2H_0},
\end{align*}
where $\tilde{\partial}_{t}=(
                 \frac{\partial}{\partial t_{1}},
\frac{1}{2}\frac{\partial}{\partial t_{2}},
\frac{1}{3}\frac{\partial}{\partial t_{3}},\cdots)$.
\end{prop}

\begin{lem}\cite{djkm1,djkm2} \label{Vi}Let $V_1(z)=\sqrt{2}\;\Phi\,  
\phi(-z)\psi^*(z)\,\Phi^{-1},\,
V_0(z)=\sqrt{2}\;\Phi\, \phi(z)\psi(z)\,\Phi^{-1}.$
Then we have
\begin{align*}
V_1(z)=e^{-\xi(t,z)}e^{2\,\xi(\widetilde{\partial_t},z^{-1})}\theta  
e^{-\alpha}z^{-2H_0},\quad
V_0(z)=e^{\xi(t,z)}e^{-2\,\xi(\widetilde{\partial_t},z^{-1})}\theta  
e^{\alpha} z^{2H_0}.
\end{align*}
\end{lem}

Due to Lemma \ref{Vi}, we can write the actions of $f_i$ on  
$\C[t]\otimes_\C\pmb{\Omega}$
in terms of  formal contour integrals
\begin{equation*}
f_0
=\sqrt{-1}^{\,-1}\oint z^{-1}V_0(z)dz,\quad
f_1
=-\oint V_1(z)dz,\label{vo}
\end{equation*}
where  we set $\oint A(z) dz=A_{-1}$ for $A(z)=\sum_{n}A_nz^n$.

\begin{lem}\label{V^l}
$$
V_i(z_\ell)\cdots V_i(z_2)V_i(z_1)
=
(-1)^{\frac{\ell(\ell-1)}{2}}\Delta(z)^2
e^{(-1)^i\sum_{j}\xi(t,z_j)}e^{2(- 
1)^{i+1}\sum_{j}\xi(\widetilde{\partial_t},z_j^{-1})}
\theta^\ell e^{(-1)^i\ell\alpha}(z_1\cdots z_\ell)^{2(-1)^iH_0}.
$$
Here $\Delta(z)=\det(z_{i}^{j-1})_{1 \leq i,j\leq \ell}$.
\end{lem}

\begin{proof}By $V_i^0(z)$, we denote the ``zero mode''
$\theta e^{(-1)^i\alpha}z^{2(-1)^iH_0}$ of $V_i(z).$
Then by using the relations $\theta e^{\pm\alpha}=-  
e^{\pm\alpha}\theta,$ and
$z_j^{\pm 2 H_0}e^{\pm \alpha}=z_j^2\cdot e^{\pm  
\alpha} z_j^{\pm 2 H_0},$
we have
$$V_i^0(z_\ell)\cdots V_i^0(z_2)V_i^0(z_1)=(-1)^{\frac{\ell(\ell-1)}{2}}
\left(\prod_{j=1}^\ell z_j^{2j-2}\right)\theta^\ell
e^{(-1)^i\ell\alpha} (z_1\cdots z_\ell)^{2(-1)^iH_0}.
$$
On the other hand, by the
standard calculus of vertex operators, we have
$$
V_i^{+}(z_2)V_i^{-}(z_1)=\left(1-\frac{z_1}{z_2}\right)^2V_i^{-}(z_1)  
V_i^{+}(z_2),
$$
where we set $V_i^{-}(z)=e^{(-1)^i\xi(t,z)}$,
$V_i^{+}(z)=e^{(-1)^{i+1}2\xi(\widetilde{\partial_t},z^{-1})}.$
Then the lemma follows immediately.
\end{proof}

For $\lambda\in \mathcal{P}$, we denote
by $$\mathrm{S}^\lambda(z)=\det(z_i^{\lambda_j+j-1})/\det(z_i^{j-1})$$
the Schur function with
respect to $z=(z_1,\ldots,z_\ell).$
We use the well-known orthogonality relation
$$
\frac{1}{(2\pi\sqrt{-1})^\ell}\int_{T^\ell}\mathrm{S}^\lambda(z)
\overline{\mathrm{S}^\mu(z)}\,|\Delta(z)|^2\,
\frac{dz_1}{z_1}\cdots \frac{dz_\ell}{z_\ell}
=\ell!\;\delta_{\lambda,\mu},
$$
where we denote by $T^\ell =
\{z=(z_j)\in \C^\ell\; ;\;|z_j|=1\}$, the $\ell$-dimensional torus.
Since $\overline{\mathrm{S}^\mu(z)}=\mathrm{S}^\mu(z^{-1})
=\mathrm{S}^\mu(z_1^{-1},\ldots,z_\ell^{-1})$ for $z\in T^\ell$,
we can rewrite this relation as
\begin{align}\label{ortho}
\oint\cdots\oint \mathrm{S}^\lambda(z){\mathrm{S}^\mu(z^{-1})}
(-1)^{\frac{\ell(\ell-1)}{2}}\Delta(z)^2(z_1\cdots  
z_\ell)^{-\ell}\,dz_1\cdots dz_\ell
=\ell!\;\delta_{\lambda,\mu}.
\end{align}
We also utilize the following form of the Cauchy identity:
\begin{align}\label{cauchy}
e^{\sum_{j=1}^\ell\xi(t,z_j)}&=\sum_{\ell(\lambda) \leq \ell}
\mathrm{S}^\lambda(z)S_{\lambda}(t).
\end{align}
We remark
\begin{align*}
e^{-\sum_{j=1}^\ell\xi(t,z_j)}&=\sum_{\ell(\lambda) \leq \ell}
(-1)^{|\lambda|}\mathrm{S}^\lambda(z)S_{\lambda'}(t).
\end{align*}
Here $\lambda'$ is the conjugate of $\lambda$.

Put
$$f_i^{(\ell)}=\frac{1}{\ell!} f_{i}^{\ell}\ (i=0,1).$$
\begin{lem}\label{rect} For $m>0$, we have
\begin{align}
f_0^{(\ell)}\theta^me^{-m\alpha}\nonumber
&=
\zeta'_{-m,\ell,1}
S_{\square(\ell,2m+1-\ell)}(t)\,\theta^{m+\ell}e^{( 
\ell-m)\alpha}\nonumber,
\end{align}
where
$$
\zeta'_{-m,\ell,0}=\sqrt{-1}^{(2m-1)\ell}\;(1\leq \ell\leq 2m-1).\quad
$$
Similarly we have
\begin{align}
f_1^{(\ell)}\theta^me^{m\alpha}\nonumber
&=
\zeta'_{m,\ell,1}
S_{\square(2m-\ell,\ell)}(t)\,\theta^{m+\ell}e^{( 
m-\ell)\alpha}\nonumber,
\end{align}
where
$$
\zeta'_{m,\ell,1}=(-1)^{m\ell}.
$$
\end{lem}

\begin{proof}
In view of the relation $e^{\ell\alpha}\cdot\theta^m=(-1)^{\ell  
m}\theta^m\cdot e^{\ell\alpha}$, we have
by Lemma \ref{V^l}
$$
V_0(z_\ell)\cdots V_0(z_1)\,\theta^me^{-m\alpha}
=(-1)^{\ell  
m}(-1)^{\frac{\ell(\ell- 
1)}{2}}\Delta(z)^2e^{\sum_j{\xi(t,z_j)}}\theta^{m+\ell}e^{(- 
m+\ell)\alpha}(z_1\cdots z_\ell)^{-2m}.
$$
Using this, we have
\begin{align}
&f_0^{(\ell)}\theta^me^{-m\alpha}\nonumber\\
&=\frac{\sqrt{-1}^{-\ell}}{\ell!}
\oint\cdots\oint (-1)^{\frac{\ell(\ell-1)}{2}}\Delta(z)^2
(z_1\cdots z_\ell)^{-2m-1}e^{\sum_j{\xi(t,z_j)}}dz_1\cdots dz_\ell
\cdot(-1)^{\ell m}\theta^{m+\ell}e^{(-m+\ell)\alpha}\nonumber\\
&=\sqrt{-1}^{-\ell}(-1)^{\ell  
m}S_{\square(\ell,2m+1-\ell)}(t)\,\theta^{m+\ell}e^{(- 
m+\ell)\alpha}\nonumber
\end{align}
where we carried out the contour integral
by using (\ref{ortho}) and (\ref{cauchy}).
Here we remark that 
$${\rm S}^{\square(\ell,m)}(z)=(z_{1}z_{2}\cdots z_{\ell})^{m}.$$

In a similar way, we have
$
\zeta'_{m,\ell,1}=(-1)^{m\ell}.
$
We just note that  
$S_{\square(2m-\ell,\ell)}(-t)=(-1)^{\ell(2m-\ell)}S_{\square(\ell,2m- 
\ell)}(t).$
Detail of the calculation is left to
the reader.
\end{proof}

The following pictures express the $f_{0}$- and $f_{1}$-action to
each maximal weight.\\

\begin{center}
\includegraphics[width=6cm]{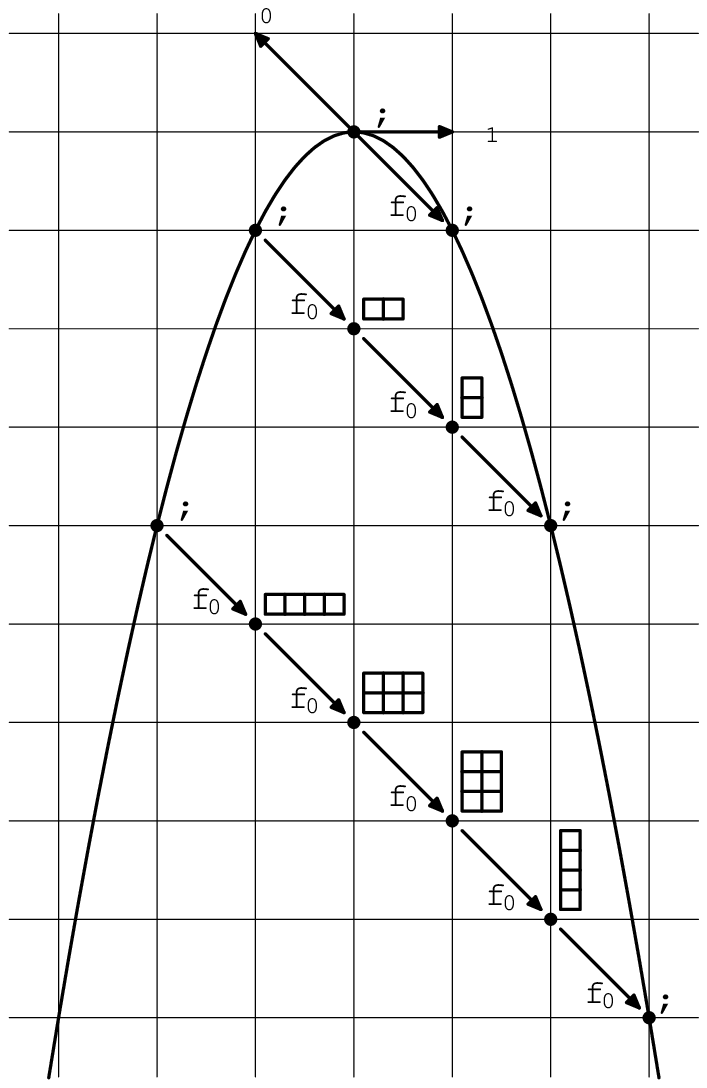}
\hspace{2cm}
\includegraphics[width=6cm]{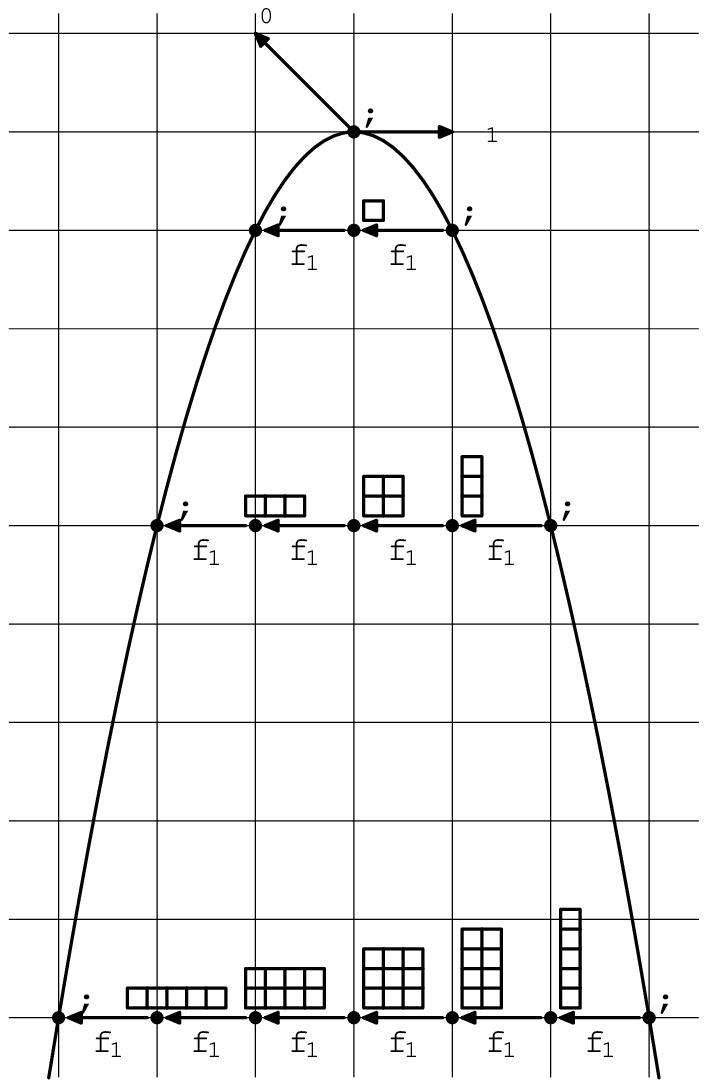}
\section{Proof of the main theorem}\label{proof}
\end{center}

First we have
\begin{align*}
\Phi\left(f_i^{(\ell)}|c_m\rangle\right)&=\sqrt{2}^{- 
\varepsilon_m}\sum_{\lambda \in {I}_i^\ell(c_m)}
\Phi\left(\sqrt{2}^{\,a(\lambda)}|\lambda\rangle\right)\quad
\nonumber\\
&=\sqrt{2}^{-\varepsilon_m}\sum_{\lambda \in  
{I}_i^\ell(c_m)}
\zeta_{m,\ell,i}(\lambda)\,Q_{\lambda[0]}(t)S_{\lambda[1]}(t^{(2)})\,\theta 
^{m+\ell}
\,e^{(m+(-1)^i\ell)\alpha}.
\end{align*}

Second we have seen in the previous section that 
     $$
     f_1^{(\ell)} \Phi ( | c_m\rangle)
      =
    \sqrt{2}^{-\varepsilon_m} \sqrt{-1}^{-m} \zeta'_{m,\ell,1}
     S_{\square(2m-\ell,\ell)}(t) \theta^{m+\ell} e^{(m-\ell)\alpha}
     $$
and
     $$
     f_0^{(\ell)} \Phi ( | c_{-m} \rangle)
      =
     \sqrt{2}^{-\varepsilon_m}\sqrt{-1}^{-m}  \zeta'_{-m,\ell,0}
     S_{\square(\ell,2m-\ell+1)}(t) \theta^{m+\ell} e^{(\ell-m)\alpha},
     $$
for $m > 0$.
Therefore all we have to show is the following:
\begin{lem}\label{signdet}
For $\lambda \in I_{i}^{\ell}(c_{m})$, we have
$$
\zeta_{m,\ell,i}(\lambda)=\sqrt{-1}^{-|m|}\zeta'_{m,\ell,i}\delta(\lambda).
$$
\end{lem}
%
%
%
%
%
We prove this lemma together with looking at some examples for help. 
We set 
$$|\lambda\rangle=\beta_{\lambda_1}\cdots\beta_{\lambda_{2s}}|\emptyset\rangle,$$
where $\lambda_1>\cdots>\lambda_{2s}\geq 0.$
If we ignore the factor $\sqrt{-1}$ in (\ref{kakikae}), the set
$\{\beta_{\lambda_1},\ldots,\beta_{\lambda_{2s}}\}$ is
decomposed into
the three parts
$$
\mathcal{I}=\{\psi_{i_1},\ldots,\psi_{i_N}\},
\quad
\mathcal{J}=\{\psi^*_{j_1},\ldots,\psi^*_{j_{N^*}}\},
\quad
\mathcal{K}=\{\phi_{k_1},\ldots,\phi_{k_a}\},
$$
where
 $a$, $N$ and $N^*$ the number of 
$\phi$'s, $\psi$'s and $\psi^*$'s respectively 
and
$i_1>\cdots>i_N\geq 0>j_1>\cdots>j_{N^*}$, 
$k_1>\cdots>k_a\geq 0.$
Actually, $I=\{i_1,\ldots,i_N\}$ ( resp. $J=\{j_1,\ldots,j_{N^*}\}$) is  
nothing but
$\iota (\pmb{\lambda}^{(1)})$ (resp. $\iota^*(\pmb{\lambda}^{(3)})$),
and $K=\{{k_1},\ldots,{k_a}\}$ corresponds to  
$\pmb{\lambda}^{(0)}\cup
\pmb{\lambda}^{(2)}.
$
According to the following operations, we shall rewrite $|\lambda\rangle$ into
 its ``normal form''
such as in Lemma \ref{QS}.

\begin{enumerate}
\item[\{OP. 0\}]
Rewrite  $\beta$'s into $\phi$'s, $\psi$'s and $\psi^*$'s according to (\ref{kakikae}). 
\item[\{OP. 1\}]
Rewrite the vacuum $|\emptyset\rangle$ into $\psi_{-1} \psi_{-2} \cdots \psi_{j_{N^{*}}}|0,j_{N^*} \rangle
$, i.e.,
$$|\emptyset\rangle=\psi_{-1} \psi_{-2} \cdots \psi_{j_{N^{*}}}|0,j_{N^*} \rangle.$$ 
\item[\{OP.  2\}]
Repeat the following operations in order of $m=j_{1},j_{2},\cdots,j_{N^{*}}$:
\begin{center}
Move $\psi^*_{m}$ to the left side of $\psi_{m}$ and
remove $\psi^*_{m}$ by using the relation
$\psi^*_{m}\psi_{m}=1-\psi_{m}\psi^*_{m}.$
\end{center}
\item[\{OP. 3\}]
Move $\phi_{j}$'s to the  left of $\psi$'s in order of $j=k_{1},k_{2}\cdots,k_{a}$.
\end{enumerate}
We divide  our argument into four cases.

\underline{Case 1}: We consider the case of $i=1$ and $m=2n>0.$
We have 
$$\zeta'_{m,\ell,1}=(-1)^{m\ell}=1.$$
Our purpose is to  rewrite $|\lambda \rangle$ into its normal form and compute a factor $\zeta_{m,\ell,i}(\lambda)$.
We employ the following example for our understanding;
$$m=6, \ell=6\ {\rm and}\ \mu=(21,19,13,10,7,2) \in I_{0}^{6}(c_{6}).$$
In this case, since all elements of $K$ are odd,    
 \{OP. 0\} gives a factor 
$\sqrt{-1}^{m}$.
In our example, we have
$$\beta_{21}\beta_{19}\beta_{13}\beta_{10}\beta_{7}\beta_{2}|\emptyset\rangle
\stackrel{\rm\{OP. 0\}}{=}\sqrt{-1}^6 \psi_{5}\psi^*_{-5}\psi_{3}\phi_{3}\psi^*_{-2}\phi_{1}|\emptyset\rangle.$$
We neglect the factor $\sqrt{-1}^m$
for the moment.
After rewriting the vacuum according to \{OP. 1\}, 
 we move $\psi^*_{j_1}$ to the
left side of $\psi_{j_{1}}$.
Then $\psi^*_{j_1}$ jumps $-j_{1}-1$ elements of $\mathcal{I}\cup\mathcal{K}$
and $\psi_{-1},\cdots,\psi_{j_{1}+1}$.
Therefore this operation gives a factor 
$$(-1)^{(-j_{1}-1)+(-j_{1}-1)}=1.$$
We apply  this operation to $\psi^*_{j_{m}}$ in oder of $m=2,3,\cdots,{N^*}$. Then we 
have a factor 
$$(-1)^{(-j_{m}-m)+(-j_{m}-m)}=1$$
for each $m$.
Therefore \{OP. 2\} gives a factor 
$$(-1)^{\sum_{m=1}^{N^*}2(-j_{m}-1)}=1.$$
In our example, we have
\begin{align*}
\sqrt{-1}^6 \psi_{5}\psi^*_{-5}\psi_{3}\phi_{3}\psi^*_{-2}\phi_{1}|\emptyset\rangle
&
\stackrel{\rm\{OP. 1\}}{=}
\sqrt{-1}^6 \psi_{5}\psi^*_{-5}\psi_{3}\phi_{3}\underline{\underline{\psi^*_{-2}}}\phi_{1}\psi_{-1}\psi_{-2}
\psi_{-3}\psi_{-4}\psi_{-5}|0,-5\rangle\\
&
\stackrel{\rm\{OP. 2\}}{=}
(-1)^{1+1}\sqrt{-1}^6 \psi_{5}\underline{\underline{\psi^*_{-5}}}\psi_{3}\phi_{3}\phi_{1}\psi_{-1}
\psi_{-3}\psi_{-4}\psi_{-5}|0,-5\rangle\\
&
\stackrel{\rm\{OP. 2\}}{=}
(-1)^{1+1}(-1)^{3+3}\sqrt{-1}^6 \psi_{5}\psi_{3}\phi_{3}\phi_{1}\psi_{-1}
\psi_{-3}\psi_{-4}|0,-5\rangle\\
&
=
\sqrt{-1}^6 \psi_{5}\psi_{3}\phi_{3}\phi_{1}\psi_{-1}
\psi_{-3}\psi_{-4}|0,-5\rangle.
\end{align*}
From the Definition \ref{gabacus},  
 the factor occurred by \{OP. 3\} is  $\delta(\lambda)$.
We see this fact through our example. We compute 
\begin{align*}
\sqrt{-1}^6 \psi_{5}\psi_{3}\underline{\underline{\phi_{3}}}\phi_{1}\psi_{-1}
\psi_{-3}\psi_{-4}|0,-5\rangle
&\stackrel{\rm\{OP. 3\}}{=}
\sqrt{-1}^6(-1)^2 \phi_{3}\psi_{5}\psi_{3}\underline{\underline{\phi_{1}}}\psi_{-1}
\psi_{-3}\psi_{-4}|0,-5\rangle\\
&\stackrel{\rm\{OP. 3\}}{=}
\sqrt{-1}^6(-1)^{2+2} \phi_{3}\phi_{1}\psi_{5}\psi_{3}\psi_{-1}
\psi_{-3}\psi_{-4}|0,-5\rangle
\end{align*}
and obtain  $\delta(\mu)=(-1)^{2+2}$ from 4-bar abacus of $\mu$;
\begin{align*}
 {\begin{array}{ccc}
0&1&3\\
\maru{2}&&\\
4&5&\maru{7}\\
6&&\\
8&9&11\\
\maru{10}&&\\
12&\maru{13}&15\\
14&&\\
18&17&\maru{19}\\
20&&\\
22&\maru{21}&23\\
\end{array}}.
\end{align*}
Now we obtain  
$$|\lambda\rangle=
\delta(\lambda)
\sqrt{-1}^{\,m}
\phi_{k_1}\cdots
\phi_{k_a}
\psi_{i_1}\cdots \psi_{i_{N}}
\psi_{-1}\psi_{-2}\cdots\widehat{\psi_{j_{1}}}
\cdots
\widehat{\psi_{j_{2}}}
\cdots
\widehat{\psi_{j_{N^*}}}
|0,j_{N^*}\rangle.$$
and $\zeta_{m,\ell,1}(\lambda)=\delta(\lambda)\sqrt{-1}^{m}$. 
Since
$m$ is even, we have $\zeta_{m,\ell,1}(\lambda)=\delta(\lambda)\sqrt{-1}^{-m}\zeta'_{m,\ell,1},$

\underline{Case 2}: We consider the case of $i=1$ and $m=2n+1>0.$ 
Then we have
$$\zeta'_{m,\ell,1}=(-1)^{m\ell}=(-1)^\ell.$$
The only difference
from the case 1  is the existence of  $\phi_0=\beta_0$  
in the
right end of $\beta$'s.
The element $\phi_{0}$ causes a factor $(-1)^{N+N^*}$,
because $\phi_{0}$ is jumped by the elements of $\mathcal{J}$ (\{OP. 2\}) and jump  the elements of $\mathcal{I}$ 
(\{OP. 3\}).
For example,
if $\mu=(19,13,10,7,2)\in I_{1}^{6}(c_{5})$, then we have
\begin{align*}
|\mu\rangle
&=
\beta_{19}\beta_{13}\beta_{10}\beta_{7}\beta_{2}\beta_{0}|\emptyset\rangle\\
&\stackrel{\rm\{OP. 0\}}{=}\sqrt{-1}^5\psi^*_{-5}\psi_{3}\phi_{3}\psi^*_{-2}\phi_{1}\phi_{0}|\emptyset\rangle\\
&\stackrel{\rm\{OP. 1\}}{=}\sqrt{-1}^5\psi^*_{-5}\psi_{3}\phi_{3}\psi^*_{-2}\phi_{1}\phi_{0}\psi_{-1}\psi_{-2}
\psi_{-3}\psi_{-4}\psi_{-5}|0,-5\rangle .
\end{align*}
In this example $\phi_{0}$ causes a factor $(-1)^{2+3}$.
Therefore we have
$$
|\lambda\rangle
=\delta(\lambda)\sqrt{-1}^{\,m}(-1)^{N+N^*}
\phi_{k_1}\cdots
\phi_{k_a}
\psi_{i_1}\cdots \psi_{i_{N}}
\psi_{-1}\psi_{-2}\cdots\widehat{\psi_{j_{1}}}
\cdots
\widehat{\psi_{j_{2}}}
\cdots
\widehat{\psi_{j_{N^*}}}
|0,j_{N^*}\rangle.
$$
By using the relations 
$$
\begin{cases}
m=N+N^*+a\\
\ell=a+2N^*,
\end{cases}$$
we have $(-1)^{N+N^*}$ is equal to $(-1)^{m-\ell}.$
Thus we have $\zeta_{m,\ell,1}(\lambda)=\delta(\lambda)\sqrt{-1}^{\,m}(-1)^{m-\ell}$
and
$\zeta_{m,\ell,i}(\lambda)=\delta(\lambda)\sqrt{-1}^{-|m|}\zeta'_{m,\ell,i},$ 
since $m$ is odd.

The following case 3 and case 4 are for $i=0$. 
In these cases, we should put  a bead on 0 of the 4-bar abacus
when $\lambda_{m+1}=0\ (m>0)$.
If we take care of this point, we can rewrite $|\lambda\rangle$
into its normal form
and 
determine $\zeta_{m,\ell,0}(\lambda)$
as same as the case 1 and case 2.   

\underline{Case 3}: 
We consider the case of $i=0$ and $m=2n+1>0.$ 
Then we have 
$$\zeta'_{-m,\ell,0}=\sqrt{-1}^{(2m-1)\ell}=\sqrt{-1}^{\ell}.$$
Take an element $\lambda\in {I}_0^\ell(c_{-m}).$ Remark that $\lambda_{m+1}=0$ or $1$.
We further divide this case into subcases:
\begin{enumerate}
\item[(a)]
 $\lambda_{m+1}=0,$\  i.e., $|c_{-m}\rangle$ has the end term $\beta_{0}=\phi_{0}$,
 \item[(b)]
  $\lambda_{m+1}=1,$\  i.e., $|c_{-m}\rangle$ does not have the end term $\beta_{0}=\phi_{0}$.
  \end{enumerate}
Let us first consider the subcase (a).
We have
$$
|\lambda\rangle=\delta(\lambda)\sqrt{-1}^{(m+1)-a}(-1)^{N^*}
\phi_{k_1}\cdots
\phi_{k_a}
\psi_{i_1}\cdots \psi_{i_{N}}
\psi_{-1}\psi_{-2}\cdots\widehat{\psi_{j_{1}}}
\cdots
\widehat{\psi_{j_{2}}}
\cdots
\widehat{\psi_{j_{N^*}}}
|0,j_{N^*}\rangle$$
and
$$
\zeta_{-m,\ell,0}(\lambda)=\delta(\lambda)\sqrt{-1}^{\,m+1-a}(-1)^{N^*}.
$$
Here we remark 
 the factor $(-1)^{N^*}$ caused by $\phi_0$ (\{OP. 2\})
and
 the factor $(-1)^{N}$ caused by $\phi_{0}$ is included by $\delta(\lambda)$ (\{OP. 3\}).
 For example, set $\mu=(12,9,3,0)\in I_{0}^{4}(c_{-3})$. Then we compute 
 \begin{align*}
 |\mu\rangle&=\beta_{12}\beta_{9}\beta_{3}\beta_{0}|\emptyset\rangle\\
 &\stackrel{\rm\{OP. 0\}}{=}
 \sqrt{-1}^{2}\phi_{6}\psi_{2}\psi^*_{-1}\phi_{0}|\emptyset\rangle\\
&\stackrel{\rm\{OP. 1\}}{=}
\sqrt{-1}^{2}\phi_{6}\psi_{2}\psi^*_{-1}\phi_{0}\psi_{-1}|0,-1\rangle\\
&\stackrel{\rm\{OP. 2\}}{=}
\sqrt{-1}^{2}(-1)^1\phi_{6}\psi_{2}\phi_{0}|0,-1\rangle\\
&\stackrel{\rm\{OP. 3\}}{=}
\sqrt{-1}^{2}(-1)^1(-1)^1\phi_{6}\phi_{0}\psi_{2}|0,-1\rangle
 \end{align*}
 and $\delta(\mu)=(-1)^{1+0}$ from the 4-bar abacus of $\mu$
 \begin{align*}
 {\begin{array}{ccc}
\maru{0}&1&\maru{3}\\
{2}&&\\
4&5&{7}\\
6&&\\
8&\maru{9}&11\\
{10}&&\\
\maru{12}&13&15
\end{array}}.
\end{align*}
Now we have the relations
$$\begin{cases}
N-N^*=-m+\ell\\
a+N+N^*=m+1.
\end{cases}$$
We eliminate $N$ from these equations to get $a+2N^*=2m-\ell+1.$
Then we have 
$$\zeta_{-m,\ell,0}(\lambda)=\delta(\lambda)\sqrt{-1}^{-m+\ell},$$ which is
equal to $\delta(\lambda)\sqrt{-1}^{-m}\zeta'_{-m,\ell,0}.$
For subcase (b), we have $\psi_0$ instead of the
absence of $\phi_0.$
So the same factor $(-1)^{N^*}$ occurs
when we exchange $\psi_0$ and $\psi^*$'s.
Then the formula is the same as (a).
The readers can check  this fact by using  an example $\mu=(12,9,3,1)\in I_{0}^{4}(c_{-3})$.

\underline{Case 4}: We consider the case of $i=0$ and $m=2n>0.$ 
Then we have
$$\zeta'_{-m,\ell,0}=\sqrt{-1}^{(2m-1)\ell}=\sqrt{-1}^{-\ell}.$$ 
Consider $\lambda\in {I}_0^\ell(c_{-m}).$
We further divide this case into the subcases:
\begin{enumerate}
 \item[(c)]
  $\lambda_{m+1}=1$ i.e., $|\lambda\rangle$ ends with  
$\psi_0\phi_0$
\item[ (d)] 
$\lambda_{m+1}=0$
i.e.,
$|\lambda\rangle$ does not contain $\phi_0$ nor $\psi_0.$
\end{enumerate}
In the case (c), 
 we remark, by the similar argument of the case 1, \{OP. 2\} does not cause any sign change 
and \{OP. 3\} causes $\delta(\lambda)(-1)^N$.
We have
$$
|\lambda\rangle=\sqrt{-1}^{(m+2)-a}\delta(\lambda)(-1)^{N}
\phi_{k_1}\cdots
\phi_{k_a}
\psi_{i_1}\cdots \psi_{i_{N}}\psi_{-1}\psi_{-2}\cdots\widehat{\psi_{j_{1}}}
\cdots
\widehat{\psi_{j_{2}}}
\cdots
\widehat{\psi_{j_{N^*}}}
|0,j_{N^*}\rangle$$
and
$$
\zeta_{-m,\ell,0}(\lambda)=\sqrt{-1}^{\,m+2-a}\delta(\lambda)(-1)^{N}.
$$
For example, if $\mu=(15,13,9,4,1) \in I_{0}^{6}(c_{-4})$, then we compute
  \begin{align*}
 |\mu\rangle&=\beta_{15}\beta_{13}\beta_{9}\beta_{4}\beta_{1}\beta_{0}|\emptyset\rangle\\
 &\stackrel{\rm\{OP. 0,1\}}{=}
 \sqrt{-1}^{4}\underline{\underline{\psi^*_{-4}}}
 \psi_{3}\psi_{2}\phi_{2}\psi_{0}\phi_{0}\psi_{-1}\psi_{-2}\psi_{-3}\psi_{-4}|0,-4\rangle\\
&\stackrel{\rm\{OP. 2\}}{=}
 \sqrt{-1}^{4}(-1)^{3+2+3}\psi_{3}\psi_{2}
 \underline{\underline{\phi_{2}}}
 \psi_{0}\phi_{0}\psi_{-1}\psi_{-2}\psi_{-3}|0,-4\rangle\\
&\stackrel{\rm\{OP. 3\}}{=}
 \sqrt{-1}^{4}(-1)^2\phi_{2}\psi_{3}\psi_{2}\psi_{0}
 \underline{\underline{\phi_{0}}}
 \psi_{-1}\psi_{-2}\psi_{-3}|0,-4\rangle\\
&\stackrel{\rm\{OP. 3\}}{=}
 \sqrt{-1}^{4}(-1)^2(-1)^3\phi_{2}\phi_{0}\psi_{3}\psi_{2}\psi_{0}\psi_{-1}\psi_{-2}\psi_{-3}|0,-4\rangle\\
 \end{align*}
 and $\delta(\mu)=(-1)^2$ from the 4-bar abacus:
 \begin{align*}
 {\begin{array}{ccc}
{0}&\maru{1}&{3}\\
{2}&&\\
\maru{4}&5&{7}\\
6&&\\
8&\maru{9}&11\\
{10}&&\\
{12}&\maru{13}&\maru{15}
\end{array}}.
\end{align*}
Since we have 
$$
\begin{cases}
a+2N=\ell+2 \\
 m\ \text{is\ odd},
 \end{cases}$$
we can see $\zeta_{-m,\ell,0}=\delta(\lambda)\sqrt{-1}^{-m}\zeta'_{-m,\ell,0}.$
The subcase (d) is the most cumbersome one.
By the definition, we include 
the sign $(-1)^N$ in $\delta(\lambda)$ arising from the exchanges of the dummy ``$\phi_0$''
and $\psi$'s.
So we have to compensate the same factor to
get $$\zeta_{-m,\ell,0}(\lambda)=\sqrt{-1}^{m-a}\delta(\lambda)(-1)^N.$$
The readers can check  this fact by using an example $\mu=(15,13,9,4)\in I_{0}^{5}(c_{-4})$.
Now using $a+2N=\ell,$ we have $\zeta_{-m,\ell,0}(\lambda)
=\delta(\lambda)\sqrt{-1}^{m-\ell}$.
 Since  $m$ is even, we have
 $\sqrt{-1}^{m-\ell}=\sqrt{-1}^{-m-\ell}=\sqrt{-1}^{-m}\zeta'_{m,\ell,0}$.

\end{document}